\documentclass[12pt]{amsart}

\usepackage{amsmath,amsfonts,amssymb,amsthm}
\usepackage[english]{babel}
\usepackage{latexsym}
\usepackage{amscd, epsfig}

\numberwithin{equation}{section}

\DeclareMathOperator{\comm}{{comm}}
\DeclareMathOperator{\cf}{{cf}}
\DeclareMathOperator{\rtq}{{rq}}

\DeclareMathOperator{\qq}{\mathbf{q}}

\newcommand{\ts}{{\hskip.04cm}}

\newcommand{\bu}{\bullet}

\newcommand{\Si}{\Sigma}

\newcommand{\ga}{\gamma}
\newcommand{\si}{\sigma}
\newcommand{\de}{\delta}

\newcommand{\al}{\alpha}
\newcommand{\be}{\beta}

\newcommand{\vp}{\varphi}

\newcommand{\wt}{\widetilde}

\newcommand{\tth}{{\hspace{.04cm}}}
\newcommand{\pP}{ \Sigma}

\newcommand{\Z}{ \mathbb Z}
\newcommand{\N}{ \mathbb N}

\newcommand{\C}{ \mathbb C}
\newcommand{\p}[1]{\mathcal{#1}}
\newcommand{\s}[1]{\mathbf{#1}}
\newcommand{\set}[1]{\{#1\}}

\newcommand{\detq}{{\textstyle \det_q}}
\newcommand{\detqq}{{\textstyle \det_{\s q}}}
\DeclareMathOperator{\sdet}{sdet}

\DeclareMathOperator{\inv}{inv}
\DeclareMathOperator{\oinv}{oinv}
\DeclareMathOperator{\cyc}{cyc}

\newtheorem{thm}{Theorem}[section]
\newtheorem{lemma}[thm]{Lemma}

\newtheorem{prop}[thm]{Proposition}

\newtheorem{exm}[thm]{Example}
\newtheorem{rem}[thm]{Remark}
\newtheorem{Warn}[thm]{Caution}

    {\end{Example}}
    {\end{Remark}}

\newtheoremstyle{prf}
  {0pt}{9pt}{}{}{\itshape}{.}{.5em}{}
\theoremstyle{prf}
\newtheorem*{prf}{Proof}

\title[Non-commutative MacMahon Master Theorem]
{Non-commutative extensions of the MacMahon Master Theorem}

\author{Matja\v z Konvalinka$^\star$ and Igor Pak$^\star$}

\thanks{\thinspace ${\hspace{-1.35ex}}^\star$Department of Mathematics,
MIT, Cambridge, MA.
\hskip.03cm
Email:
\texttt{\{konvalinka,pak\}@math.mit.edu}}

\begin{document}

\begin{abstract}
We present several non-commutative extensions of the MacMahon
Master Theorem, further extending the results of Cartier-Foata
and Garoufalidis-L\^e-Zeilberger.  The proofs are combinatorial
and new even in the classical cases.  We also give applications
to the $\beta$-extension and Krattenthaler-Schlosser's
$q$-analogue.
\end{abstract}

\maketitle

\section*{Introduction} \label{s:intro}
The MacMahon Master Theorem is one of the jewels in enumerative
combinatorics,  and it is as famous and useful as it is
mysterious.  Most recently, a new type of algebraic generalization
was proposed in~\cite{garoufalidis} and was further
studied in \cite{foata1,foata2,foata3,hai}.  In this paper
we present further generalizations of the MacMahon Master Theorem
and several other related results.  While our generalizations are
algebraic in statement, the heart of our proofs is completely
bijective, unifying all generalizations.  In fact, we give
a new bijective proof of the (usual) MacMahon Master Theorem,
modulo some elementary linear algebra.  Our approach seems to
be robust enough to allow further generalizations in this
direction.


Let us begin with a brief outline of the history of the subject.
The Master Theorem was discovered in 1915 by Percy MacMahon in
his landmark two-volume \emph{``Combinatory Analysis''}, where
he called it ``a Master Theorem in the Theory of
Partitions''~\cite[page 98]{macmahon}.  Much later, in the early
sixties, the real power of Master Theorem was discovered,
especially as a simple tool for proving binomial
identities (see~\cite{goulden}).  The proof of the
MacMahon Master Theorem using Lagrange
inversion is now standard, and the result is often viewed
in the analytic context~\cite{good,goulden}.

An algebraic approach to MacMahon Master Theorem goes back to
Foata's thesis~\cite{foata-th}, parts of which were later
expanded in~\cite{cf} (see also~\cite{la}).
The idea was to view the theorem as a result on ``words''
over a (partially commutative) alphabet, so one can
prove it and generalize it by means of simple combinatorial
and algebraic considerations.  This approach became highly
influential and led to a number of new related results
(see e.g.~\cite{kob,minoux,viennot,zeilberger}).

While the Master Theorem continued to be extended in several
directions (see~\cite{foata4,krattenthaler}), the ``right''
q- and non-commutative analogues of the results evaded discovery
until recently.  This was in sharp contrast with the Lagrange
inversion, whose $q$- and non-commutative analogues were understood
fairly well~\cite{garsia,gr,gessel,gs,Kr,pak,Si}.
Unfortunately, no reasonable
generalizations of the Master Theorem followed from these results.

An important breakthrough was made by Garoufalidis, L\^e and
Zeilberger (GLZ), who introduced a new type of $q$-analogue, with a puzzling
algebraic statement and a technical proof~\cite{garoufalidis}.
In a series of papers, Foata and Han first modified and extended the
Cartier-Foata combinatorial approach to work in this algebraic setting,
obtaining a new (involutive) proof of the GLZ-theorem~\cite{foata1}.
Then they developed a beautiful ``1 = q''
principle which gives perhaps the most elegant explanation of the
results~\cite{foata2}. They also analyze a number of specializations
in~\cite{foata3}.
Most recently, Hai and Lorenz gave an interesting algebraic proof of the
GLZ-theorem, opening yet another direction for exploration
(see Section~\ref{s:final}).

\medskip

This paper presents a number of generalizations of the MacMahon Master
Theorem in the style of Cartier-Foata and Garoufalidis-L\^e-Zeilberger.
Our approach is bijective and is new even in the classical cases, where
it is easier to understand.
This is reflected in the structure of the paper: we present
generalizations one by one, gradually moving from well known results
to new ones. The paper is largely self-contained and no
background is assumed.

We begin with basic definitions, notations and statements of the
main results in Section~\ref{s:basic}.  The proof of the (usual)
MacMahon Master Theorem is given in Section~\ref{macm}.
While the proof here is elementary, it is the basis for our approach.
A straightforward extension to the Cartier-Foata case is given
in Section~\ref{cf}.  The right-quantum case is presented in
Section~\ref{rq}.  This is a special case of the
GLZ-theorem, when $q=1$.  Then we give a $q$-analogue of the
Cartier-Foata case (Section~\ref{cfq}), and the GLZ-theorem
(Section~\ref{q}).  The subsequent results
are our own and can be summarized as follows:

\smallskip

$\bu$ \ The Cartier-Foata $(q_{ij})\tth$-analogue (Section~\ref{cfqij}).

\smallskip

$\bu$ \ The right-quantum $(q_{ij})\tth$-analogue (Section~\ref{qij}).

\smallskip

$\bu$ \ The super-analogue (Section~\ref{super}).

\smallskip

$\bu$ \ The $\beta$-extension (Section~\ref{beta}).

\smallskip
\noindent
The $(q_{ij})\tth$-analogues are our main result; one of them specializes to the
GLZ-theorem when all $q_{ij} = q$.  The super-analogue is a direct
extension of the classical MacMahon Master Theorem to commuting and
anti-commuting variables.  Having been overlooked in previous
investigations, it is a special case of the $(q_{ij})\tth$-analogue,
with some $q_{ij}=1$ and others~$= -1$.  Our final extension is somewhat
tangential to the main direction, but is similar in philosophy.
We show that our proof of the MacMahon Master Theorem can be easily
modified to give a non-commutative generalization of the so called
$\be$-extension, due to Foata and Zeilberger~\cite{foata4}.

In Section~\ref{ks} we present one additional observation on the subject.
In~\cite{krattenthaler}, Krattenthaler and Schlosser obtained an
intriguing $q$-analogue of the MacMahon Master Theorem, a result which
on the surface does not seem to fit the above scheme.  We prove that in
fact it follows from the classical Cartier-Foata generalization.

As the reader shall see, an important technical part of our proof
is converting the results we obtain into traditional form.  This is
basic linear algebra in the classical case, but in non-commutative
cases the corresponding determinant identities are either less known
or new.  For the sake of completeness, we present concise proofs of
all of them in Section~\ref{proofs}.
We conclude the paper with final remarks and open problems.

\section{Basic definitions, notations and main results} \label{s:basic}

\subsection{Classical Master Theorem} \label{ss:basic-words}
We begin by stating the Master Theorem in the classical form:

\begin{thm}  \label{macm1}
{\rm (MacMahon Master Theorem)} \
 Let $A=(a_{ij})_{m \times m}$ be a complex matrix,
 and let $x_1,\ldots,x_m$ be a set of variables.
 Denote by $G(k_1,\ldots,k_m)$ the coefficient of
 $x_1^{k_1}\cdots x_m^{k_m}$ in
 \begin{equation} \label{e:macm-G}
 \prod_{i=1}^m (a_{i1}x_1+\ldots+a_{im}x_m)^{k_i}.
\end{equation}
 Let $t_1,\ldots,t_m$ be another set of variables, and
 $T = (\de_{ij} t_i)_{m \times m}$.
 Then
 \begin{equation} \label{macm2}
  \sum_{(k_1,\ldots,k_m)} \, G(k_1,\ldots,k_m)
  \ t_1^{k_1} \cdots t_m^{k_m} \, =
  \, \frac 1{\det(I - TA)},
 \end{equation}
 where the summation is over all nonnegative integer
 vectors $(k_1,\ldots,k_m)$.
\end{thm}

By taking $t_1=\ldots=t_m=1$ we get
\begin{equation} \label{e:macm-D}
\sum_{(k_1,\ldots,k_m)} \, G(k_1,\ldots,k_m) \, =
  \, \frac 1{\det(I - A)}\,,
\end{equation}
whenever both sides of the equation are well defined, for example
when all $a_{ij}$ are formal variables.  Moreover,
replacing $a_{ij}$ in~\eqref{e:macm-D} with $a_{ij}\,t_i$
shows that~\eqref{e:macm-D} is actually equivalent
to~\eqref{macm2}.  We will use this observation throughout
the paper.

\subsection{Non-commuting variables} \label{ss:basic-noncom}

Consider the following algebraic setting.  Denote
by~$\p A$ the algebra (over~$\C$) of formal power series with
non-commuting variables $a_{ij}$, $1 \leq i,j \leq m$.
Elements of $\p A$ are infinite linear combinations of words
in variables $a_{ij}$ (with coefficients in~$\C$).
In most cases we will take
elements of $\p A$ modulo some ideal~$\p I$ generated by
a finite number of relations.  For example, if~$\p I$ is
generated by $a_{ij}a_{kl} = a_{kl}a_{ij}$ for all
$i,j,k,l$, then $\p A/\p I$ is the symmetric algebra
(the free commutative algebra with $m^2$ variables $a_{ij}$,
$1 \leq i,j \leq m$).

\medskip

Throughout the paper we assume that $x_1,\ldots,x_m$ commute
with~$a_{ij}$, and that $x_i$ and $x_j$ commute up to
some nonzero complex weight, i.e.\hspace{-0.07cm} that
$$
 x_j x_i \, = \, q_{ij} \,x_i x_j\,, \ \ \mbox{for all} \ \ i < j
$$
with $q_{ij} \in \C$, $q_{ij} \ne 0$.
We can then expand the expression
\begin{equation} \label{e:macm-G-arrow}
\prod_{i =\tth 1..{m}}^{\longrightarrow}
\bigl(a_{i1}x_1+\ldots+a_{im}x_m\bigr)^{k_i},
\end{equation}
move all $x_i$'s to the right and order them.  Along the way,
we will exchange pairs of variables $x_i$ and $x_j$, producing
a product of~$q_{ij}$'s.  We can then extract
the coefficient at $x_1^{k_1}\cdots x_m^{k_m}$.  As before,
 we will denote this coefficient by $G(k_1,\ldots,k_m)$.
 Each such coefficient will be a finite sum of products
 of a monomial in $q_{ij}$'s, $1 \leq i < j \leq m$,
 and a word $a_{i_1j_1} \, \ldots \,a_{i_{\ell}j_\ell}$,
such that $i_1 \le \ldots \le i_\ell$, the number of
variables $a_{i,\ast}$ is equal to~$k_i$, and the number of
variables $a_{\ast,j}$ is equal to~$k_j$.

\medskip

To make sense of the right-hand side of \eqref{e:macm-D} in the
non-commutative case, we need to generalize the determinant.
Throughout the paper the (non-commutative) determinant
will be given by the formula
\begin{equation} \label{e:macm-Det}
\det(B) \, = \, \sum_{\si \in S_m} \, w(\si) \, b_{\si_1{1}} \cdots \ts
b_{\si_m m} \,,
\end{equation}
where $\si = (\si_1,\ldots,\si_m)$ is a permutation and $w(\si)$ is
a certain constant weight of~$\si$.  Of course,
$w(\si) = (-1)^{\inv(\sigma)}$ is the usual case, where $\inv(\sigma)$
is the number of inversions in~$\sigma$.

Now, in all cases we consider the weight of the identity permutation
will
be equal to $1$: $w(1,\ldots,m) = 1$.  Substituting $B = I-A$
in~\eqref{e:macm-Det}, this gives us
$$\frac1{\det(I-A)} \, = \, \frac1{1 - \Si} \, = \, 1 \, + \, \Si \,
+ \, \Si^2 \, + \, \ldots\,,
$$
where~$\Si$ is a certain finite sum of words in~$a_{ij}$ and both left and
right inverse of $\det(I-A)$ are equal to the infinite sum on the right.
From now on, whenever justified, we will always use the fraction notation
as above in non-commutative situations.

\medskip

In summary, we just showed that both
$$\sum_{(k_1,\ldots,k_m)} \, G(k_1,\ldots,k_m) \ \ \ \text{and} \ \ \,
\frac 1{\det(I - A)}$$
are well-defined elements of $\p A$.  The generalizations of the
Master Theorem we present in this paper will state that these
two expressions are equal modulo a certain ideal~$\p I$.
In the classical case, the MacMahon Master Theorem gives
that for the ideal
$\p I_{\comm}$ generated by $a_{ij}a_{kl} = a_{kl}a_{ij}$, for all
$1\le i,j,k,l \le m$.

\subsection{Main theorem} \label{ss:basic-main}
Fix complex numbers $q_{ij} \ne 0$, where $1 \le i < j \le m$.
Suppose the variables $x_1,\ldots,x_m$ are $\qq$-commuting:
\begin{equation} \label{e:qij-x}
 x_j \ts x_i \, = \, q_{ij} \,x_i \ts x_j\,, \ \ \text{for all} \ \ i < j,
\end{equation}
and that they commute with all $a_{ij}$.  Suppose also that
the variables $a_{ij}$ $\qq$-commute within columns:
\begin{equation} \label{e:qij-a1}
a_{jk}\ts a_{ik} \, = \, q_{ij} \, a_{ik}\ts a_{jk}\,, \ \
\text{for all} \ \ i < j,
\end{equation}
and in addition satisfy the following quadratic equations:
\begin{equation} \label{e:qij-a2}
a_{jk}\, a_{il} \, - \, q_{ij} \, a_{ik}\ts a_{jl} \, + \,
q_{kl} \, a_{jl}\ts a_{ik} \, - \,
q_{kl} \ts q_{ij} \, a_{il}\ts a_{jk} \, = 0 \,, \ \
\text{for all} \ \ i < j, \ k <l.
\end{equation}
We call $A = (a_{ij})$ with entries satisfying~\eqref{e:qij-a1}
and~\eqref{e:qij-a2} a \emph{right-quantum $\s q$-matrix}.

For a matrix $B = (b_{ij})_{m \times m}$, define the
$\s q$-\emph{determinant} by
\begin{equation} \label{e:qq-det}
\detqq (b_{ij})=\sum_{\sigma} \, w(\si) \,
 b_{\sigma_1 1}\cdots b_{\sigma_m m},
\end{equation}
where
$$w(\si) = \prod_{i<j, \ \sigma_i>\sigma_j}
(-q_{\sigma_j\sigma_i})^{-1}.
$$

\begin{thm} \label{t:main}
Let $A =(a_{ij})_{m \times m}$ be a right-quantum $\s q$-matrix.
Denote the
coefficient of $x_1^{k_1}\cdots x_m^{k_m}$ in
 $$\prod_{i =\tth 1..{m}}^{\longrightarrow} \,
 \bigl(a_{i1}x_1+\ldots+a_{im}x_m\bigr)^{k_i}.$$
 by $G(k_1,\ldots,k_m)$. Then
 \begin{equation} \label{e:main}
  \sum_{(k_1,\ldots,k_m)} \, G(k_1,\ldots,k_m)=\frac 1{\detqq(I-A)},
 \end{equation}
where the summation is over all nonnegative integer
 vectors $(k_1,\ldots,k_m)$.
\end{thm}

Theorem~\ref{t:main} is the ultimate extension of the classical
MacMahon Master Theorem.  Our proof of the theorem uses a
number of technical improvements which become apparent in
special cases.  While the proof is given in Section~\ref{qij},
it is based on all previous sections.

\section{A combinatorial proof of the MacMahon Master Theorem}
\label{macm}

\subsection{Determinant as a product}
\label{ss:master-setup}
Let $B = (b_{ij})$ be an invertible $m \times m$ matrix over~$\C$.
Denote by $B^{11}$ the matrix $B$ without the first row and the first
column, by $B^{12,12}$ the matrix $B$ without the first two rows and
the first two columns, etc.  For the entries of the inverse matrix we
have:
\begin{equation} \label{e:macm-frac}
\bigl(B^{-1}\bigr)_{11} \, = \frac{\det B^{11}}{\det B}.
\end{equation}
Substituting $B = I - A$ and iterating~\eqref{e:macm-frac}, we obtain:
$$\aligned
& \left(\frac{1}{I-A}\right)_{11}
  \left(\frac{1}{I-A^{11}}\right)_{22}
  \left(\frac{1}{I-A^{12,12}}\right)_{33} \,
\cdots \,\ts \frac{1}{1-a_{mm}} \\
& \quad = \, \frac{\det\left(I-A^{11}\right)}{\det(I-A)}\ts \cdot\ts
\frac{\det\left(I-A^{12,12}\right)}{\det(I-A^{11})} \ts \cdot \ts
\frac{\det\left(I-A^{123,123}\right)}{\det(I-A^{12,12})}\, \cdots\ts
\frac{1}{1-a_{mm}}\\
& \quad =\, \frac 1{\det(I-A)}\,,
\endaligned
$$
provided that all minors are invertible.  Now let
$a_{ij}$ be commuting variables as in
Subsection~\ref{ss:basic-words}.  We obtain that the
right-hand side of equation~\eqref{e:macm-D} is the product
of entries in the inverses of matrices, and we need to prove
the following identity:
\begin{equation} \label{e:macm-D1}
\sum \,
G(k_1,\ldots,k_m)  =
 \left(\frac{1}{I-A}\right)_{11}
  \left(\frac{1}{I-A^{11}}\right)_{22}
\left(\frac{1}{I-A^{12,12}}\right)_{33}
\cdots \, \frac{1}{1-a_{mm}}\,.
\end{equation}

Since $(I-A)^{-1} = I + A + A^2 + \ldots$, we get a combinatorial
interpretation of the $(11)$-entry:
\begin{equation} \label{e:macm-11}
 \left(\frac{1}{I-A}\right)_{11} \, = \, \sum \, a_{1j_1} a_{j_1j_2}
 \cdots a_{j_{\ell}1}\,,
\end{equation}
where the summation is over all finite sequences
$(j_1,\ldots,j_\ell)$, where $j_r \in \{1,\ldots,m\}$,
$1 \le r \le \ell$.  A combinatorial interpretation of
the other product terms is analogous.  Recall that we already have
a combinatorial interpretation of $G(k_1,\ldots,k_m)$ as a
summation of words.  Therefore, we have reduced
the Master Theorem to an equality between two summations
of words~\eqref{e:macm-D}, where all the summands have
a positive sign. 
To finish the proof we construct
an explicit bijection between the families of words corresponding to
both sides.

\subsection{The bijection} \label{ss:master-paths}
Throughout the paper we consider \emph{lattice steps}
of the form $(x,i) \to (x+1,j)$ for some $x,i,j \in \Z$,
$1 \leq i,j \leq m$.  We think of $x$ being drawn along $x$-axis,
increasing from left to right, and refer to~$i$ and~$j$ as
the \emph{starting height} and \emph{ending height}, respectively.

From here on, we represent the step $(x,i) \to (x+1,j)$ by the
variable~$a_{ij}$. Similarly, we represent a finite sequence of
steps by a word in the alphabet $\set{a_{ij}}$, $1 \leq i,j \leq m$,
i.e.\hspace{-0.07cm} by an element of algebra $\p A$.
If each step in a sequence
starts at the ending point of the previous step, we call such
a sequence a \emph{lattice path}.
\smallskip

Define a \emph{balanced sequence} (\emph{b-sequence}) to
be a finite sequence of steps
\begin{equation} \label{e:star} \al \, = \, \bigl\{
(0,i_1) \to (1,j_1) \,, \, (1,i_2) \to (2,j_2)\,, \,\ldots \,,
\, (\ell-1,i_{\ell}) \to (\ell,j_\ell)\bigr\}\ts,
\end{equation}
such that the number of steps starting at height~$i$ is equal
to the number of steps ending at height~$i$, for all~$i$.
We denote this number by $k_i$, and call $(k_1,\ldots,k_m)$
the \emph{type} of the b-sequence.
Clearly, the total number of steps in the path $\ell = k_1+\ldots+k_m$.

\smallskip

Define an \emph{ordered sequence} (\emph{o-sequence})
to be a b-sequence where the steps starting at smaller height
always precede steps starting at larger heights. In other words,
an o-sequence of type $(k_1,\ldots,k_m)$ is a sequence of $k_1$
steps starting at height~$1$, then $k_2$ steps starting at
height~$2$, etc., so that $k_i$ steps end at height $i$.
Denote by $\s O(k_1,\ldots,k_m)$ the set of all o-sequences
of type $(k_1,\ldots,k_m)$.

\smallskip

Now consider a lattice
path from $(0,1)$ to $(x_1,1)$ that never goes below $y=1$ or
above $y=m$, then a lattice path from $(x_1,2)$ to $(x_2,2)$
that never goes below $y=2$ or above $y=m$, etc.; in the end,
take a straight path from $(x_{m-1},m)$ to $(x_m,m)$.  We
will call this a \emph{path sequence} (\emph{p-sequence}).
Observe that every p-sequence is also a b-sequence.
Denote by
$\s P(k_1,\ldots,k_m)$ the set of all p-sequences of type
$(k_1,\ldots,k_m)$.

\begin{figure}[ht!]
\begin{center}
 \includegraphics{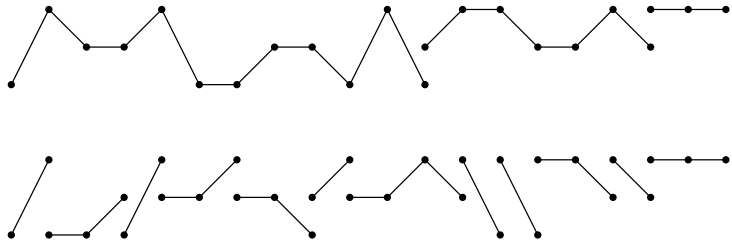}
 \caption{An o-sequence and a p-sequence of type $(4,7,8)$}
 \label{fig1}
\end{center}
\end{figure}

\begin{exm} \, {\rm
Figure \ref{fig1} presents the o-sequence
associated with the word
$$
a_{13}a_{11}a_{12}a_{13}a_{22}a_{23}a_{22}a_{21}a_{23}
a_{22}a_{23}a_{32}a_{31}a_{31}a_{33}a_{32}a_{32}a_{33}a_{33}
$$
and the p-sequence associated with
$$a_{13}a_{32}a_{22}a_{23}a_{31}a_{11}a_{12}a_{22}a_{21}
a_{13}a_{31}a_{23}a_{33}a_{32}a_{22}a_{23}a_{32}a_{33}a_{33}.
$$
}
\end{exm}

We are ready now to establish a connection between balanced sequences
and the equation~\eqref{e:macm-D1}.  First, observe that choosing
a term of
$$\prod_{i=1}^m (a_{i1}x_1+\ldots+a_{im}x_m)^{k_i}
$$
means choosing a term $a_{1*}x_*$ $k_1$ times, then choosing
a term $a_{2*}x_*$ $k_2$ times, etc., and then multiplying
all these terms.  In other words, each term on the left-hand side
of~\eqref{e:macm-D1} corresponds to an o-sequence in
$\s O(k_1,\ldots,k_m)$ for a unique vector $(k_1,\ldots,k_m)$.
Similarly, by~\eqref{e:macm-11}, a term on the right-hand side
of~\eqref{e:macm-D1} corresponds to a p-sequence,
i.e.\hspace{-0.07cm} to an element
of $\s P(k_1,\ldots,k_m)$ for a unique vector $(k_1,\ldots,k_m)$.

Let us define a bijection
$$\vp \, : \, \s O(k_1,\ldots,k_m) \, \longrightarrow \, \s P(k_1,\ldots,k_m)
$$
with the property that the word $\vp(\alpha)$ is a rearrangement
of the word~$\alpha$, for every o-sequence $\alpha$.

Take an o-sequence $\alpha$, and let $[0,x]$ be the maximal
interval on which it is part of a p-sequence, i.e.\hspace{-0.07cm} the
maximal interval $[0,x]$ on which the o-sequence has the property
that if a step ends at level $i$, and the following step starts at
level $j > i$, the o-sequence stays on or above height $j$ afterwards.
Let~$i$ be the height at~$x$.
Choose the step $(x',i) \to (x'+1,i')$ in the o-sequence that is
the first to the right of $x$ that starts at level $i$
(such a step exists because an o-sequence is a balanced sequence).
Continue switching this step with the one to the left until it
becomes the step $(x,i) \to (x+1,i')$.  The new object is part of a
p-sequence at least on the interval $[0,x+1]$.  Continuing this
procedure we get a p-sequence $\vp(\alpha)$.

For example, for the o-sequence given in Figure~1 we have $x=1$
and~$i=3$.  The step we choose then is $(12,3) \to (13,1)$,
i.e.\hspace{-0.07cm} $x' = 12$.

\begin{lemma}\label{l:vp}
 The map~$\vp : \s O(k_1,\ldots,k_m) \to \s P(k_1,\ldots,k_m)$
 constructed above is a bijection.
\end{lemma}
\begin{prf}
 Since the above procedure never switches two steps that begin at the
 same height, there is exactly one o-sequence that maps into a given
 b-sequence: take all steps starting at height $1$ in the b-sequence
 in the order they appear, then all the steps starting at height $2$
 in the p-sequence in the order they appear, etc. Clearly, this map
 preserves the type of a b-sequence. \qed
\end{prf}

\begin{figure}[ht!]
\begin{center}
 \includegraphics{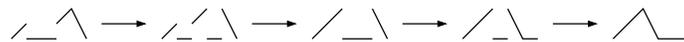}
 \caption{Transforming an o-sequence into a p-sequence.}
 \label{fig3}
\end{center}
\end{figure}

\begin{exm} \, {\rm
Figure~\ref{fig3}
shows the switches for an o-sequence of type $(3,1,1)$, and the p-sequence
in Figure~\ref{fig1} is the result of applying this procedure to
the o-sequence in the same figure (we need $33$ switches). }
\end{exm}

In summary, Lemma~\ref{l:vp} establishes the desired bijection between two sides
of equation~\eqref{e:macm-D1}.  This completes the proof
of the theorem. \qed

\subsection{Refining the bijection}  Although we already established
the MacMahon Master Theorem, in the next two subsections we will
refine and then elaborate on the proof.  This will be useful
when we consider various generalizations and modifications of
the theorem.

First, let us define \emph{q-sequences} to be the b-sequences we get
in the transformation of an o-sequence into a p-sequence with the
above procedure (including the o-sequence and the p-sequence).
Examples of q-sequences can be seen in Figure~\ref{fig3}, where
an o-sequence is transformed into a p-sequence via the intermediate
q-sequences.

Formally, a q-sequence is a b-sequence with the following properties:
it is part of a p-sequence on some interval $[0,x]$ (and this part ends at some
height~$i$); the rest of the sequence has non-decreasing starting heights,
with the exception of the first step to the right of~$x$ that starts
at height~$i$, which can come before some steps starting at lower levels.
For a q-sequence~$\alpha$, denote by $\psi(\alpha)$ the
q-sequence we get by performing the switch defined above; for a
p-sequence $\alpha$ (where no more switches are needed),
$\psi(\alpha)=\alpha$.  By construction, map~$\psi$ always switches
steps that start on different heights.

For a balanced sequence \eqref{e:star},
define the \emph{rank}~$r$ as follows:
$$
r\, := \, \bigl|\{(s,t)~:~i_s >i_t, \, 1\le s < t \le \ell \}\bigr|\,.
$$
Clearly, o-sequences are exactly the balanced sequences of rank~0.
Note also that the map~$\psi$ defined above increases by~$1$ the rank
of sequences that are not p-sequences.

Write $\s Q_n(k_1,\ldots,k_m)$ for the union of two sets of
b-sequences of type $(k_1,\ldots,k_m)$: the set of
all q-sequences with rank $n$ and the set of p-sequences with
rank $< n$; in
particular, $\s O(k_1,\ldots,k_m)=\s Q_0(k_1,\ldots,k_m)$ and
$\s P(k_1,\ldots,k_m)=\s Q_N(k_1,\ldots,k_m)$ for $N$ large
enough (say, $N \geq \binom{\ell}{2}$ will work).

\begin{lemma}
 The map $\psi: \s Q_n(k_1,\ldots,k_m) \to \s Q_{n+1}(k_1,\ldots,k_m)$
 is a bijection for all~$n$.
\end{lemma}
\begin{prf}
 A q-sequence of rank $n$ which is not a p-sequence is mapped into
 a q-sequence of rank $n+1$, and $\psi$ is the identity map on
 p-sequences. This proves that $\psi$ is indeed a map from
 $\s Q_n(k_1,\ldots,k_m)$ to $\s Q_{n+1}(k_1,\ldots,k_m)$.
 It is easy to see that $\psi$ is injective and surjective.\qed
\end{prf}

The lemma gives another proof that
$\vp = \psi^N : \s O(k_1,\ldots,k_m) \to \s P(k_1,\ldots,k_m)$
is a bijection.  This is the crucial observation which will be
used repeatedly in the later sections.

\medskip

Let us emphasize the importance of bijections~$\psi$ and~$\vp$
in the language
of ideals.  Obviously we have $\psi(\alpha)=\alpha$ modulo $\p I_{\comm}$
for every q-sequence $\alpha$.  Consequently, $\vp(\alpha)=\alpha$
modulo $\p I_{\comm}$ for every o-sequence, and we have
$$\sum \ \vp(\alpha)\, =\, \sum \ \alpha \  \mod \ \p I_{\comm}\,,
$$
where the sum is over all o-sequences~$\al$.
From above, this can be viewed as a restatement of the
MacMahon Master Theorem~\ref{macm1}.

\subsection{Meditation on the proof} \
The proof we presented above splits into two (unequal) parts:
combinatorial and linear algebraic.  The combinatorial part
(the construction of the bijection~$\vp$) is the heart of the proof
and will give analogues of~\eqref{e:macm-D1} in non-commutative cases
as well.  While it is fair to view the equation~\eqref{e:macm-D1} as
the ``right'' generalization of the Master Theorem, it is preferable
if the right-hand side is the inverse of some version of the determinant, for
both aesthetic and traditional reasons.  This is also how our
Main Theorem~\ref{t:main} is stated.

The linear algebraic part, essentially the
equation~\eqref{e:macm-frac}, is trivial
in the commutative (classical) case. The non-commutative
analogues we consider are much less trivial, but largely known.
In the most general case considered in the Main Theorem
the formula follows easily from the results of Manin on quantum
determinants~\cite{manin2,manin-book} and advanced technical
results of Etingof and Retakh who proved~\eqref{e:macm-frac}
for quantum determinants~\cite{ER} in a more general setting
(see further details in Section~\ref{s:final}).

To avoid referring the technicalities to other people's work and
deriving these basic linear algebra facts from much more general
results, we include our own proofs of the analogues
of~\eqref{e:macm-frac}. These proofs are moved to Section~\ref{proofs} and
we try to keep them as concise and elementary as possible.

\section{The Cartier-Foata case} \label{cf}

In this section, we will assume that the variables $x_1,\ldots,x_m$
commute with each other and with all $a_{ij}$, and that
\begin{equation} \label{cf4}
 a_{ij}\ts a_{kl}\,=\,a_{kl}\ts a_{ij} \ \ \text{for all} \ \, i \neq k\ts.
\end{equation}
The matrix $A = (a_{ij})$ which satisfies the
conditions above is called a \emph{Cartier-Foata matrix}.

For any matrix $B=(b_{ij})_{m \times m}$ (with non-commutative entries)
define the \emph{Cartier-Foata determinant}:
$$\det B\, =\,\sum_{\sigma\in S_m}
(-1)^{\inv(\sigma)} \, b_{\sigma_1 1}\cdots
 b_{\sigma_m m}\,.
$$
Note that the order of terms in the product is important in general,
though not for a Cartier-Foata matrix.

\begin{thm}[Cartier-Foata] \label{cf1}
Let $A = (a_{ij})_{m \times m}$ be a Cartier-Foata matrix.
 Denote by $G(k_1,\ldots,k_r)$ the coefficient of
 $x_1^{k_1}\cdots x_m^{k_m}$ in the product
 $$\prod_{i =\tth 1..{m}}^{\longrightarrow} \,
 \bigl(a_{i1}x_1+\ldots+a_{im}x_m\bigr)^{k_i}.
 $$
 Then
 \begin{equation} \label{cf2}
  \sum_{(k_1,\ldots,k_m)}  G(k_1,\ldots,k_m) \, = \,
  \frac 1{\det(I-A)}\,,
 \end{equation}
 where the summation is over all nonnegative integer
 vectors $(k_1,\ldots,k_m)$, and $\det(\cdot)$ is the
 Cartier-Foata determinant.
 \end{thm}

Clearly, Theorem~\ref{cf1} is a generalization of the MacMahon
Master Theorem~\ref{macm1}.  Let us show that our proof of the
Master Theorem easily extends to this case.
We start with the following well known
technical result (see e.g.~\cite{foata5}).

\begin{prop} \label{cf3}
 If $A=(a_{ij})_{m \times m}$ is a Cartier-Foata matrix, then
  $$
  \left(\frac1{I-A}\right)_{11} \, = \, \frac1{\det(I-A)} \, \cdot \,
  \det\left(I-A^{11}\right),
  $$
where $\det(\cdot)$ is the Cartier-Foata determinant.
\end{prop}
For completeness, we include a straightforward proof of the
proposition in Section~\ref{proofs}.

\begin{proof}[Proof of Theorem~\ref{cf1}]
Denote by $\p I_{\cf}$ the ideal generated by relations
$a_{ij}a_{kl}=a_{kl}a_{ij}$ for all $1\le i,j,k,l \le m$,
with~$i \neq k$.
Observe that the terms of the left-hand side of~\eqref{cf2} correspond
to o-sequences.  Similarly, by Proposition~\ref{cf3} and
equation~\eqref{e:macm-11}, the terms on the right-hand side correspond
to p-sequences. Therefore, to prove the theorem it suffices to
show that
\begin{equation} \label{e:cf-I}
\sum \,\alpha \,= \, \sum \,\vp(\alpha) \  \mod \  \p I_{\cf} \,,
\end{equation}
where the sum is over all o-sequences of a fixed type $(k_1,\ldots,k_m)$.

As mentioned earlier, all switches we used in the
construction of~$\psi$
involve steps starting at different heights. This means that for a
q-sequence~$\alpha$, we have
$$\psi(\alpha)\,=\,\alpha \ \mod \  \p I_{\cf}\,,
$$
which implies~\eqref{e:cf-I}.  This completes the proof of
the theorem.  \end{proof}

\section{The right-quantum case} \label{rq}

In this section, we will assume that the variables $x_1,\ldots,x_m$
commute with each other and with all $a_{ij}$, and that we have
\begin{eqnarray}
 a_{jk}\ts a_{ik} \, & = & \, a_{ik}\ts a_{jk}, \label{rq1} \\
 a_{ik}\ts a_{jl}\,-\,a_{jk}\ts a_{il} \, & =
 & \, a_{jl}\ts a_{ik}\, - \, a_{il}\ts a_{jk}, \label{rq2}
\end{eqnarray}
for all $1 \le i,j,k,l \le m$. We call $A=(a_{ij})_{m \times m}$ whose
entries satisfy these relations a \emph{right-quantum matrix}.

Note that a Cartier-Foata matrix is a right-quantum matrix.
The following result is an important special case of the
GLZ-theorem (Theorem~\ref{q}) and a generalization of
Theorem~\ref{cf1}.

\begin{thm} \label{rq4}
Let $A = (a_{ij})_{m\times m}$ be a right-quantum matrix.
 Denote by $G(k_1,\ldots,k_r)$ the coefficient of
 $x_1^{k_1}\cdots x_m^{k_m}$ in the product
 $$
 \prod_{i =\tth 1..{m}}^{\longrightarrow} \,
 \bigl(a_{i1}x_1+\ldots+a_{im}x_m\bigr)^{k_i}.
 $$
 Then
 \begin{equation} \label{rq5}
  \sum_{(k_1,\ldots,k_m)}  G(k_1,\ldots,k_m) \, =
  \, \frac 1{\det(I-A)}\,,
 \end{equation}
 where the summation is over all nonnegative integer
 vectors $(k_1,\ldots,k_m)$, and $\det(\cdot)$ is the
 Cartier-Foata determinant.
\end{thm}

Let us show that our proof of the Master Theorem extends
to this case as well, with some minor modifications.
We start with the following technical result generalizing
Proposition~\ref{cf3}.

\begin{prop} \label{rq3}
If $A = (a_{ij})$ is a right-quantum matrix, then
 $$
  \left(\frac1{I-A}\right)_{11} \, =
  \, \frac1{\det(I-A)} \, \cdot \,
  \det\left(I-A^{11}\right).
  $$
\end{prop}
For completeness, we include a proof of the proposition in
Section~\ref{proofs}.

\begin{proof}[Proof of Theorem~\ref{rq4}]
Denote by~$\p I_{\rtq}$ the ideal of~$\p A$
generated by the relations~\eqref{rq1} and~\eqref{rq2}.
As before, the proposition implies that the right-hand side of~\eqref{rq5}
enumerates all p-sequences, and it is again obvious that the left-hand side
of~\eqref{rq5}
enumerates all o-sequences. Note that it is no longer true that for
an o-sequence $\alpha$, $\vp(\alpha)=\alpha$ modulo~$\p I_{\rtq}$.
However, it suffices to prove that
\begin{equation}\label{e:rq-I}
\sum \,\vp(\alpha) \, = \, \sum \, \alpha \  \mod \ \p I_{\rtq},
\end{equation}
where the sum goes over all o-sequences
$\al \in \s O(k_1,\ldots,k_m)$.
We show this by making switches in the construction
of~$\vp$ \emph{simultaneously}.

Take a q-sequence~$\alpha$. If $\alpha$ is a p-sequence, then
$\psi(\alpha)=\alpha$. Otherwise, assume that $(x-1,i) \to (x,k)$ and
$(x,j) \to (x+1,l)$ are the steps to be switched in order to get
$\psi(\alpha)$. If $k = l$, then $\psi(\alpha)=\alpha$ modulo $\p I_{\rtq}$
by \eqref{rq1}. Otherwise, denote by $\beta$ the sequence we get by
replacing these two steps with $(x-1,i) \to (x,l)$ and $(x,j) \to (x+1,k)$.
The crucial observation is that $\beta$ is also a q-sequence, and that
its rank is equal to the rank of $\alpha$. Furthermore,
$\alpha + \beta = \psi(\alpha)+\psi(\beta) \mod \p I_{\rtq}$ because
of~\eqref{rq2}. This implies that $\sum \psi(\alpha)=\sum \alpha  \mod \p I_{\rtq}$ with
the sum over all sequences in $\s Q_n(k_1,\ldots,k_m)$.
From here we obtain~\eqref{e:rq-I} and conclude the proof of
the theorem. \end{proof}

\begin{figure}[ht!]
\begin{center}
 \includegraphics{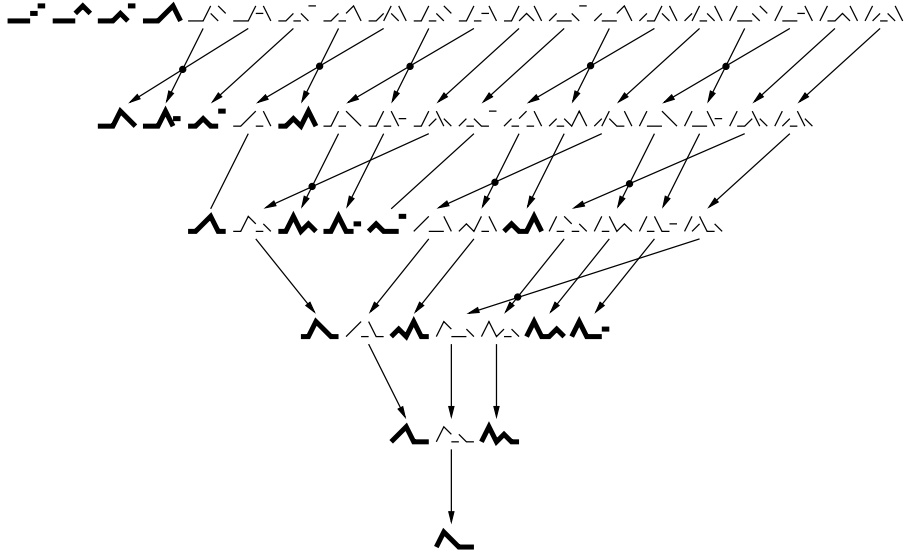}
 \caption{Transforming o-sequences into p-sequences via a series of
 simultaneous switches.}
 \label{fig2}
\end{center}
\end{figure}

\begin{exm} {\rm
Figure \ref{fig2} provides a graphical illustration for $k_1=3$,
$k_2=1$, $k_3=1$; here p-sequences are drawn in bold, an arrow from a
q-sequence $\alpha$ of rank $n$ to a q-sequence of rank $n+1$ $\alpha'$
means that $\alpha'=\psi(\alpha)$ and $\alpha'=\alpha \mod \p I_{\rtq}$,
and arrows from q-sequences $\alpha,\beta$ of rank $n$ to q-sequences
$\alpha',\beta'$ of rank $n+1$ whose intersection is marked
by a dot mean that $\alpha'=\psi(\alpha)$, $\beta'=\psi(\beta)$,
and $\alpha' + \beta' = \alpha + \beta  \mod \p I_{\rtq}$. }
\end{exm}

\section{The Cartier-Foata $q$-case} \label{cfq}

In this section, we assume that variables $x_1,\ldots,x_m$ satisfy
\begin{equation} \label{cfq1}
 x_j \ts x_i \, = \, q \, x_i \ts x_j\, \ \ \text{for}  \ \ i < j,
\end{equation}
where $q \in \C$, $q \ne 0$, is a fixed complex number.
Suppose also that $x_1,\ldots,x_m$  they commute with all~$a_{ij}$
and that we have:
\begin{eqnarray}
 a_{jl}\ts a_{ik} \, &  = & \, a_{ik}\ts a_{jl}\, \ \
 \text{for}  \ \ i < j, \, k < l\ts,
 \label{cfq2} \\
 a_{jl}\ts a_{ik} \, & = & \, q^2 \ts a_{ik} \ts a_{jl}\,,
 \ \ \text{for}  \ \ i < j, \, k > l\ts, \label{cfq3} \\
 a_{jk}\ts a_{ik} \, & = & \, q \ts a_{ik}\ts a_{jk}\,,
 \ \ \text{for}  \ \ i < j\ts.
  \label{cfq4}
\end{eqnarray}
Let us call such a matrix
$A = (a_{ij})$ a \emph{Cartier-Foata $q$-matrix}.
As the name suggests, when $q=1$ the Cartier-Foata $q$-matrix becomes
a Cartier-Foata matrix.

For a matrix $B=(b_{ij})_{m \times m}$ with non-commutative entries,
define a \emph{quantum determinant} (\emph{q-determinant}) by the
following formula:
 $$
 \detq B \, = \, \sum_{\sigma \in S_m} \, (-q)^{-\inv(\sigma)} \,
 b_{\sigma_1 1} \cdots b_{\sigma_m m}
 $$
The following result is another
important special case of the GLZ-theorem and a generalization of the
Cartier-Foata Theorem~\ref{cf1}.

\begin{thm} \label{cfq5}
Let $A = (a_{ij})_{m\times m}$ be a Cartier-Foata $q$-matrix.
 Denote by $G(k_1,\ldots,k_r)$
 the coefficient of $x_1^{k_1}\cdots x_m^{k_m}$ in
 $$\prod_{i =\tth 1..{m}}^{\longrightarrow} \,
  \bigl(a_{i1}x_1+\ldots+a_{im}x_m\bigr)^{k_i}.$$
 Then
 \begin{equation} \label{cfq6}
  \sum_{(k_1,\ldots,k_m)} G(k_1,\ldots,k_m) \,
   = \, \frac 1{\detq(I-A)}\,,
 \end{equation}
 where the summation is over all nonnegative integer
 vectors $(k_1,\ldots,k_m)$.
\end{thm}

The proof of the theorem is a weighted analogue of the proof of
Theorem~\ref{cf1}.  The main technical difference is essentially
bookkeeping of the powers of~$q$ which appear after switching
the letters~$a_{ij}$ (equivalently, the
lattice steps in the $q$-sequences).  We begin with some helpful
notation which will be used throughout the remainder of the paper.

We abbreviate the product $a_{\lambda_1\mu_1}\cdots a_{\lambda_n\mu_n}$
to $a_{\lambda,\mu}$ for $\lambda=\lambda_1\cdots \lambda_n$ and
$\mu=\mu_1\cdots \mu_n$, where~$\lambda$ and~$\mu$ are regarded as words
in the alphabet $\set{1,\ldots,m}$. For any such word $\nu=\nu_1 \cdots \nu_n$,
define the \emph{set of inversions}
$$\p I(\nu) \, = \, \set{(i,j) \colon i < j, \nu_i > \nu_j},$$
and let~$\inv\nu = |\p I(\nu)|$.  

\begin{figure}[ht!]
\begin{center}
 \includegraphics{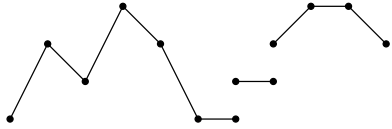}
 \caption{A p-sequence with weight $q^6$.}
\label{fig4}
\end{center}
\end{figure}

\begin{proof}[Proof of Theorem~\ref{cfq5}]
Denote by $\p I_{q-\cf}$ the ideal of~$\p A$ generated by
relations~\eqref{cfq2} -- \eqref{cfq4}.
When we expand the product
$$\prod_{i =\tth 1..{m}}^{\longrightarrow} \,
  \bigl(a_{i1}x_1+\ldots+a_{im}x_m\bigr)^{k_i},$$
move the $x_i$'s to the right and order them, the coefficient at
$a_{\lambda,\mu}$
is $q^{\inv \mu}$. This means that $\sum G(k_1,\ldots,k_m)$ is a
weighted sum of o-sequences, with an o-sequence $a_{\lambda,\mu}$ weighted
by $q^{\inv \mu}=q^{\inv \mu - \inv \lambda}$.

Choose a q-sequence $\alpha=a_{\lambda,\mu}$ and let
$\psi(\alpha)=a_{\lambda',\mu'}$.
Assume that the switch we perform is between steps
$(x-1,i) \to (x,k)$ and $(x,j) \to (x+1,l)$; write
$\lambda=\lambda_1ij\lambda_2$, $\mu=\mu_1kl\mu_2$,
$\lambda'=\lambda_1ji\lambda_2$, $\mu'=\mu_1lk\mu_2$.
If $i < j$ and $k < l$, we have $\inv \lambda'=\inv \lambda + 1$,
$\inv \mu'=\inv \mu + 1$.  By~\eqref{cfq2}, $\psi(\alpha)=\alpha$ modulo $\p I_{q-\cf}$ and
\begin{equation}\label{e:cfq-I}
q^{\inv \mu'-\inv \lambda'} \psi(\alpha)=q^{\inv \mu-\inv \lambda} \alpha
\  \mod \  \p I_{q-\cf}\,.
\end{equation}

Similarly, if $i < j$ and $k > l$, we have $\inv \lambda'=\inv \lambda + 1$,
$\inv \mu'=\inv \mu - 1$.  By~\eqref{cfq3}, we have $\psi(\alpha)=q^2\alpha$
modulo $\p I_{q-\cf}$, which implies equation~\eqref{e:cfq-I}.
If $i < j$ and $k = l$, we have $\inv \lambda'=\inv \lambda + 1$,
$\inv \mu'=\inv \mu$.  By \eqref{cfq4}, we have $\psi(\alpha)=q\alpha$
modulo $\p I_{q-\cf}$, which implies~\eqref{e:cfq-I} again.
Other cases are analogous.

Iterating equation~\eqref{e:cfq-I},
we conclude that if $\alpha=a_{\lambda,\mu}$ is an o-sequence and
$\vp(\alpha)=a_{\lambda',\mu'}$ is the corresponding p-sequence, then
$$
q^{\inv \mu'-\inv \lambda'} \vp(\alpha)\, =\, q^{\inv \mu-\inv \lambda} \, \alpha
\  \mod \ \p I_{q-\cf}\,.$$
Therefore,
\begin{equation}\label{e:cfq-II}
\sum_{(k_1,\ldots,k_m)} G(k_1,\ldots,k_m) \, = \,
\sum \, q^{\inv \mu-\inv \lambda} \, \alpha \  \mod \ \p I_{q-\cf}\,,
\end{equation}
where the sum on the right-hand side goes over all p-sequences
$\alpha=a_{\lambda,\mu}$.

Let us call a p-sequence \emph{primitive} if it starts at some height~$y$
and stays strictly above $y$ until the last step (when it returns to~$y$).
For example, the p-sequence in Figure~\ref{fig4} is a product of four
primitive p-sequences. For a primitive p-sequence $a_{\lambda,\mu}$ of
length~$\ell$, $\inv \mu - \inv \lambda = \ell - 1$, and for an arbitrary
p-sequence $a_{\lambda,\mu}$ of length $\ell$ that decomposes into~$n$
primitive p-sequences, $\inv \mu - \inv \lambda = \ell - n$.

Consider a matrix
\begin{equation} \label{atq} \wt A \, = \,
\begin{pmatrix}
a_{11} & a_{12} & \cdots & a_{1m} \\
qa_{21} & qa_{22} & \cdots & qa_{2m} \\
qa_{31} & qa_{32} & \cdots & qa_{3m} \\
\vdots & \vdots & \ddots & \vdots
\\ qa_{m1} & qa_{m2} & \cdots & qa_{mm}
\end{pmatrix}.
\end{equation}

Clearly $(\wt A^\ell)_{11}$ enumerates paths starting and ending at
height~$1$ weighted by $q^{\ell - n}$, where~$n$ is the number of steps
starting at height~$1$.  At this point we need the following generalization
of Proposition~\ref{cf3}.

\begin{prop} \label{cfq7}
 If $A=(a_{ij})_{m \times m}$ is a Cartier-Foata $q$-matrix,
 then
 $$\left(\frac1{I-\wt A}\right)_{11} \, = \,
 \frac1{\detq(I-A)} \, \cdot \, \detq\left(I-A^{11}\right).
 $$
\end{prop}

The proposition implies that the right-hand side of~\eqref{cfq6} in the theorem
enumerates all p-sequences, with $\alpha=a_{\lambda,\mu}$ weighted by
$q^{\inv \mu - \inv \lambda}$.  The equation~\eqref{e:cfq-II} above
shows that this is also the left-hand side of~\eqref{cfq6}.
This completes the proof of the theorem.
\end{proof}

\begin{exm} \ {\rm
For the p-sequence
$$\al \, = \, a_{13}a_{32}a_{24}a_{43}a_{31}a_{11}a_{22}a_{34}a_{44}a_{43}$$
shown in Figure~\ref{fig4}, we have
$$\inv(1324312344)=0+3+1+4+2+0+0+0+0+0=10$$
and
$$\inv(3243112443)=4+2+5+3+0+0+0+1+1+0=16.$$
Therefore, the p-sequence~$\al$ is weighted by~$q^6$.}
\end{exm}

\section{The right-quantum $q$-case} \label{q}

As in the previous section, we assume that variables
$x_1,\ldots,x_m$ satisfy
\begin{equation} \label{q3}
 x_j \ts x_i \, = \, q \ts x_i \ts x_j\, \ \ \text{for}  \ \ i < j,
\end{equation}
where $q \in \C$, $q \ne 0$ is a fixed complex number.
Suppose also that $x_1,\ldots,x_m$ commute with all~$a_{ij}$
and that we have:
\begin{eqnarray}
 a_{jk}\ts a_{ik} & = & q \ts a_{ik}a_{jk}\, \ \
 \text{for all}  \ \ i < j\ts,
 \label{q1} \\
 a_{ik}\ts a_{jl}\, -\, q^{-1} a_{jk}\ts a_{il} \, & = & \,
 a_{jl}\ts a_{ik}\, -q \ts a_{il}\ts a_{jk}\,
    \ \ \text{for all}  \ \ i < j,k<l\ts. \label{q2}
\end{eqnarray}
We call such matrix $A = (a_{ij})$ the \emph{right quantum q-matrix}.
It is easy to see that when $q=1$ we get a right quantum matrix defined
in Section~\ref{rq}.  In a different direction, every Cartier-Foata $q$-matrix
is also a right quantum q-matrix.  The following result of Garoufalidis, L\^e
and Zeilberger~\cite{garoufalidis} generalizes Theorems~\ref{rq4} and~\ref{cfq5}.

\begin{thm} [GLZ-theorem] \label{q4}
Let $A = (a_{ij})_{m\times m}$ be a right quantum q-matrix.
 Denote by $G(k_1,\ldots,k_r)$ the coefficient of
 $x_1^{k_1}\cdots x_m^{k_m}$ in
 $$\prod_{i =\tth 1..{m}}^{\longrightarrow} \,
  \bigl(a_{i1}x_1+\ldots+a_{im}x_m\bigr)^{k_i}.
  $$
 Then
 \begin{equation} \label{q5}
  \sum_{(k_1,\ldots,k_m)} G(k_1,\ldots,k_m) \, = \,
  \frac 1{\detq(I-A)}\,,
 \end{equation}
 where the summation is over all nonnegative integer
 vectors $(k_1,\ldots,k_m)$.
\end{thm}

The proof of the theorem is almost identical
to the one given in the previous section, with some modifications
similar to those in the proof of Theorem~\ref{rq4}.  

\begin{proof}[Proof of Theorem~\ref{q4}]
Denote by $\p I_{q-\rtq}$ the ideal of~$\p A$
generated by relations~\eqref{q1} and~\eqref{q2}.
Now, when we expand the product
$$
\prod_{i =\tth 1..{m}}^{\longrightarrow} \,
\bigl(a_{i1}x_1+\ldots+a_{im}x_m\bigr)^{k_i},$$
move the $x_i$'s to the right and order them, the coefficient at
$a_{\lambda,\mu}$ is $q^{\inv \mu}$.  Therefore, $\sum G(k_1,\ldots,k_m)$
is a weighted sum of o-sequences, with an o-sequence $a_{\lambda,\mu}$
weighted by $q^{\inv \mu}=q^{\inv \mu - \inv \lambda}$. Similar arguments
as before, now using \eqref{q1} and \eqref{q2} instead of
\eqref{cfq2} -- \eqref{cfq4}, show that
\begin{equation} \label{eq7}
\sum_{(k_1,\ldots,k_m)} G(k_1,\ldots,k_m) \, = \,
\sum q^{\inv \mu-\inv \lambda} \, a_{\lambda,\mu}
\  \mod \ \p I_{q-\rtq}\,,
\end{equation}
where the sum on the right-hand side is over all p-sequences $\alpha=a_{\lambda,\mu}$.
The following proposition generalizes Propositions~\ref{rq3} and~\ref{cfq7}.

\begin{prop} \label{q6}
If $A=(a_{ij})_{m \times m}$ is a right quantum q-matrix,
 then
 $$\left(\frac1{I-\wt A}\right)_{11} \, = \,
 \frac1{\detq(I-A)} \, \cdot \, \detq\left(I-A^{11}\right),
 $$
where~$\wt A$ is defined by~\eqref{atq}.
\end{prop}

The proposition is proved in Section~\ref{proofs}.
Now Theorem~\ref{q4}
follows from the proposition and equation~\eqref{eq7}.
\end{proof}

\section{The Cartier-Foata $q_{ij}$-case} \label{cfqij}

We can extend the results of the previous sections to the multiparameter
case.  Assume that variables $x_1,\ldots,x_m$ satisfy
\begin{equation} \label{cfqij1}
 x_j \ts x_i\, =\, q_{ij} \ts x_i \ts x_j\, \ \ \text{for}  \ \ i < j\ts,
\end{equation}
where $q_{ij} \in \C$, $q_{ij} \ne 0$ are fixed complex numbers,
$1 \le i < j \le m$.
Suppose also that $x_1,\ldots,x_m$  commute with all~$a_{ij}$
and that we have:
\begin{eqnarray}
 q_{kl} \, a_{jl}\ts a_{ik} \, & = & \, q_{ij} \, a_{ik}\ts a_{jl}\,
 \ \ \text{for} \ \ i < j,k < l\ts, \label{cfqij2} \\
 a_{jl}\ts a_{ik} \, & = & \, q_{ij} \ts q_{lk} \, a_{ik}\ts a_{jl}\,
 \ \ \text{for}
 \ \ i < j,k > l\ts,
 \label{cfqij3} \\
 a_{jk}\ts a_{ik} \,& = & \,q_{ij} \, a_{ik}\ts a_{jk}\,,
 \ \ \text{for}  \ \ i < j\ts.
 \label{cfqij4}
\end{eqnarray}
We call $A=(a_{ij})_{m \times m}$ whose
entries satisfy these relations a  \emph{Cartier-Foata $\s q$-matrix}.
When all $q_{ij}=q$ we obtain a Cartier-Foata $q$-matrix.
Thus the following result is a generalization of Theorem~\ref{cfq5}
and is a corollary of our Main Theorem~\ref{t:main}.

\begin{thm} \label{cfqij6}
Assume that $A = (a_{ij})_{m\times m}$ is a Cartier-Foata $\s q$-matrix.
 Denote by $G(k_1,\ldots,k_r)$ the coefficient of
 $x_1^{k_1}\cdots x_m^{k_m}$ in
 $$\prod_{i =\tth 1..{m}}^{\longrightarrow} \,
  \bigl(a_{i1}x_1+\ldots+a_{im}x_m\bigr)^{k_i}.$$
 Then
 \begin{equation} \label{cfqij7}
  \sum_{(k_1,\ldots,k_m)} G(k_1,\ldots,k_m) \, = \, \frac 1{\detqq(I-A)}\,,
 \end{equation}
 where the summation is over all nonnegative integer
 vectors $(k_1,\ldots,k_m)$ and $\detqq(\cdot)$ is the $\s q$-determinant
 defined by \eqref{e:qq-det}.
\end{thm}

\begin{rem} \label{r:cfqq} \,
{\rm
 If we define $q_{ii} = 1$ and $q_{ji} = q_{ij}^{-1}$ for $i < j$,
 we can write the conditions \eqref{cfqij2} -- \eqref{cfqij4} more
 concisely as
 \begin{equation} \label{cfqij5}
  q_{kl}\, a_{jl}\ts a_{ik} \, = \, q_{ij} \, a_{ik}\ts a_{jl},
 \end{equation}
 for all $i,j,k,l$, and $i \neq j$.

Let us note also that the definition of $\s q$-determinant $\detqq(B)$
for the minors of~$B$ has to be adapted as follows.  The weights
$q_{ij}$ always correspond to indices~$i,j$ of the entries~$b_{ij}$,
not the column and row numbers.  For example,
 $$
 \detqq \begin{pmatrix}
 a_{22} & a_{23} \\
 a_{32} & a_{33}
 \end{pmatrix}
 \, = \, a_{22}\ts a_{33} \,- \,q_{23}^{-1} \ts a_{32} \ts a_{23}. \eqno
 $$}
\end{rem}

We can repeat the proof of Theorem \ref{cfq5} almost verbatim.
This only requires a more careful ``bookkeeping'' as we need to keep
track of the set of inversions, not just its cardinality (the number
of inversions).

\medskip
\noindent
\begin{proof}[Proof of Theorem~\ref{cfqij6}]
Denote by $\p I_{\s q-\cf}$ the ideal of $\p A$ generated by
the relations \eqref{cfqij2} -- \eqref{cfqij4}.
When we expand the product
$$
\prod_{i =\tth 1..{m}}^{\longrightarrow} \,
\bigl(a_{i1}x_1+\ldots+a_{im}x_m\bigr)^{k_i},
$$
move the $x_i$'s to the right and order them, the coefficient
at $a_{\lambda,\mu}$ is equal to
$$\prod_{(i,j) \in \p I(\mu)} q_{\mu_j\mu_i}\,.
$$
This means that $\sum G(k_1,\ldots,k_m)$ is a weighted sum of o-sequences,
with an o-sequence $a_{\lambda,\mu}$ weighted by
$$\prod_{(i,j) \in \p I(\mu)} q_{\mu_j\mu_i} =
\prod_{(i,j) \in \p I(\mu)} q_{\mu_j\mu_i}
\prod_{(i,j) \in \p I(\lambda)} q_{\lambda_j\lambda_i}^{-1}.
$$
Now, the equation \eqref{cfqij5} implies that for every o-sequence
$\alpha=a_{\lambda,\mu}$ and $\vp(\alpha)=a_{\lambda',\mu'}$, we have
$$
\left(\prod_{(i,j) \in \p I(\mu')} q_{\mu'_j\mu'_i}
\prod_{(i,j) \in \p I(\lambda')} q_{\lambda'_j\lambda'_i}^{-1}\right)
\vp(\alpha)=\left(\prod_{(i,j) \in \p I(\mu)} q_{\mu_j\mu_i}
\prod_{(i,j) \in \p I(\lambda)} q_{\lambda_j\lambda_i}^{-1}\right) \alpha
\  \mod \ \p I_{\s q-\cf}.
$$

On the other hand, for a primitive p-sequence $a_{\lambda,\mu}$ starting
and ending at~$1$ we have:
$$\prod_{(i,j) \in I(\mu)} q_{\mu_j \mu_i}
\prod_{(i,j) \in I(\lambda)} q_{\lambda_j \lambda_i}^{-1}
= q_{1\mu_1}q_{1\mu_2}\cdots q_{1\mu_{n-1}}.
$$
This shows that all weighted p-sequences starting and ending at~$1$ are
enumerated by
\begin{equation} \label{e:wtA}
\left(\frac1{I-\wt A}\right)_{11}\,, \ \ \text{where} \ \ \
\wt A=\begin{pmatrix}
a_{11} & a_{12} & \cdots & a_{1m} \\
q_{12}a_{21} & q_{12}a_{22} & \cdots & q_{12}a_{2m} \\
q_{13}a_{31} & q_{13}a_{32} & \cdots & q_{13}a_{3m} \\
\vdots & \vdots & \ddots & \vdots \\
q_{1m}a_{m1} & q_{1m}a_{m2} & \cdots & q_{1m}a_{mm}
\end{pmatrix}.
\end{equation}
We need the following generalization of Proposition~\ref{cfq7}.

\begin{prop} \label{cfqij8}
If $A=(a_{ij})_{m \times m}$ is a Cartier-Foata $\s q$-matrix,
 then
 $$\left(\frac1{I-\wt A}\right)_{11} \, = \,
 \frac1{\detqq(I-A)} \, \cdot \, \detqq\left(I-A^{11}\right).
 $$
\end{prop}

\noindent
The proposition is proved in Section~\ref{proofs}.
From here, by the same logic as in the proofs above we obtain the
result.
\end{proof}

\section{The right-quantum $q_{ij}$-case (proof of Main Theorem~\ref{t:main})}
\label{qij}

First, by taking $q_{ii}=1$ and $q_{ji}=q_{ij}^{-1}$
for $j < i$, we can assume~\eqref{e:qij-x} holds for all $1 \le i,j \le m$.
Now equations~\eqref{e:qij-a1} and~\eqref{e:qij-a2} can be more succinctly
written as
\begin{equation} \label{qij4}
 a_{ik}\ts a_{jl} \, - q_{ij}^{-1} \,  a_{jk}\ts a_{il} \, = \,
 q_{kl}\bigl(q_{ij}^{-1} \, a_{jl}\ts a_{ik}\ts - \ts a_{il}\ts a_{jk}\bigr)
\end{equation}
for all $i,j,k,l$, such that~$i \neq j$.  Note that in this form
equation~\eqref{qij4} is a direct generalization of
equation~\eqref{q2} on one hand, with the $q_{ij}$'s arranged as in
equations~\eqref{cfqij2}--\eqref{cfqij4} on the other hand.

We also need the following (straightforward) generalization of
Propositions~\ref{q6} and~\ref{cfqij8}.

\begin{prop} \label{qij7}
If $A=(a_{ij})_{m \times m}$ is a right-quantum $\s q$-matrix,
 then
 $$\left(\frac1{I-\wt A}\right)_{11} \, = \,
 \frac1{\detqq(I-A)} \, \cdot \, \detqq\left(I-A^{11}\right),
 $$
where~$\wt A$ is given in~\eqref{e:wtA}.
\end{prop}

The proof of the proposition is in Section~\ref{proofs}.
From here, the proof of the Main Theorem follows verbatim the
proof of Theorem~\ref{cfqij6}.  We omit the details. \qed

\section{The super-case} \label{super}

In this section we present an especially interesting corollary of
Theorem~\ref{cfqij6}.

Fix a vector $\gamma = (\ga_1,\ldots,\ga_m) \in \Z_2^m$ and write
$\hat \imath$ for $\gamma_i$. If $\hat \imath = 0$, index $i$ is called
\emph{even}, otherwise it is called \emph{odd}.
We will assume that the variables $x_1,\ldots,x_m$ satisfy
\begin{equation} \label{s0}
x_j \ts x_i \, = \, (-1)^{\hat \imath \hat \jmath} \, x_i \ts x_j \ \
\text{for} \ \, i \ne j
\end{equation}
In other words, variables~$x_i$ and~$x_j$ commute unless they
are both odd: $\ga_i = \ga_j = 1$, in which case they anti-commute.
As before, suppose $x_1,\ldots,x_m$ commute with all~$a_{ij}$'s,
and that we have
\begin{eqnarray}
 a_{ik}\ts a_{jk} \, & = & \, (-1)^{\hat \imath \hat \jmath} \,
 a_{jk}\ts a_{ik}, \ \ \quad
 \text{for all} \ \ i \ne j,
 \label{s1} \\
 a_{ik}\ts a_{jl} \, & =  & \,
 (-1)^{\hat \imath \hat \jmath + \hat k \hat l} \,
 a_{jl}\ts a_{ik}, \ \ \text{for all} \ \ i \ne j, \ k \ne l\ts .
 \label{s2}
\end{eqnarray}
We call $A = (a_{ij})$ as above a \emph{Cartier-Foata super-matrix}.
Clearly, when $\gamma=(0,\ldots,0)$, we get the (usual) Cartier-Foata
matrix (see Section~\ref{cf}).

\medskip

For a permutation $\sigma$ of $\set{1,\ldots,m}$, we will denote
by $\oinv(\sigma)$ the number of \emph{odd inversions},
i.e.\hspace{-0.07cm} the number of pairs $(i,j)$ with
$\hat \imath = \hat \jmath = 1$, $i<j$, $\pi(i) > \pi(j)$.
For a matrix $B=(b_{ij})_{m \times m}$ define its
\emph{super-determinant} as
$$\sdet B=
 \sum_{\sigma \in S_m} \, (-1)^{\inv(\sigma) - \oinv(\sigma)} \,
 b_{\sigma_1 1}\cdots b_{\sigma_m m}.
$$

\begin{thm} \label{super3} {\rm (Super Master Theorem)}
 Let $A =(a_{ij})_{m\times m}$ be a Cartier-Foata super-matrix,
 and let $x_1,\ldots,x_m$ be as above.
 Denote by $G(k_1,\ldots,k_r)$ the coefficient of
 $x_1^{k_1}\cdots x_m^{k_m}$ in
 $$\prod_{i =\tth 1..{m}}^{\longrightarrow} \,
 \bigl(a_{i1}x_1+\ldots+a_{im}x_m\bigr)^{k_i}.$$
 Then
 \begin{equation} \label{super4}
  \sum_{(k_1,\ldots,k_m)} \, G(k_1,\ldots,k_m)=
  \frac{1}{\sdet(I-A)}\,,
 \end{equation}
 where the summation is over all nonnegative integer
 vectors $(k_1,\ldots,k_m)$.
\end{thm}

\begin{prf}
This is a special case of Theorem~\ref{cfqij6} for
$q_{ij}=(-1)^{\hat \imath \hat \jmath}$.  It is easy to
see that $\detqq(B) = \sdet(B)$ for all~$B$, by definition
of even and odd inversions.  The rest is a straightforward
verification.  \qed
\end{prf}

In conclusion, let us note that when $\gamma=(1,\ldots,1)$
we get a Cartier-Foata $q$-matrix with $q=-1$. Interestingly, here~$\sdet$
becomes a permanent.

\section{The $\beta$-extension} \label{beta}
In this section we first present an extension of MacMahon
Master Theorem due to Foata and Zeilberger, and then show
how to generalize it to a non-commutative setting.

First, assume that $a_{ij}$ are commutative variables
and let $\beta \in \N$ be a non-negative integer.
For $\s k=(k_1,\ldots,k_m)$, let $\pP (\s k)$ denote
the set of all permutations of the set
$$\set{(1,1),\ldots,(1,k_1),(2,1),\ldots,(2,k_2),
\ldots,(m,1),\ldots,(m,k_m)}\ts.
$$
For a permutation $\pi \in \pP(\s k)$, we define $\pi_{ij}:=i'$
whenever $\pi(i,j)=(i',j')$.  Define the \emph{weight} $v(\pi)$
by a word
$$
v(\pi)=\prod_{i=1..m}^{\longrightarrow} \,
\prod_{j=1..k_i}^{\longrightarrow} \, \, a_{i,\pi_{ij}}
$$
and the \emph{$\beta$-weight} $v_\be(\pi)$ by a product
$$
v_\beta(\pi)=\beta^{\cyc \pi}v(\pi)\ts,
$$
where $\cyc \pi$ is the number of cycles of the permutation~$\pi$.
For example, if
$$
\pi = \begin{pmatrix}
(1,1) & (1,2) & (1,3) & (2,1) & (3,1) \\
(2,1) & (1,2) & (1,1) & (3,1) & (1,3)
\end{pmatrix} \, \in \Si(3,1,1)\ts,
$$
then $v(\pi)=a_{12}a_{11}a_{11}a_{23}a_{31}$ and
$v_\beta(\pi)=\beta^2 \ts a_{12}a_{11}a_{11}a_{23}a_{31}$.

\medskip

By definition, the word $v(\pi)$ is always an o-sequence of type
$(k_1,\ldots,k_m)$.  Note now that the word
$\al \in \s O(k_1,\ldots,k_m)$ does not determine the
permutation~$\pi$ uniquely, since the second coordinate~$j'$
in $(i',j') = \pi(i,j)$ can take any value between~$1$
and~$k_{i'}$.  From here it follows that there are exactly
$k_1! \cdots k_m!$ permutations $\pi \in \pP(\s k)$
corresponding to a given o-sequence $\al \in \s O(k_1,\ldots,k_m)$.

Now, the (usual) MacMahon Master Theorem can be restated as
$$\text{(MMT)} \qquad
\sum_{\s k = (k_1,\ldots,k_m)} \,\frac{1}{k_1! \cdots k_m!} \,
\sum_{\pi \in \pP(\s k)} \, v(\pi)
\, = \, \frac1{\det(I-A)}\,,
$$
where the summation is over all non-negative integer vectors
$\s k = (k_1,\ldots,k_m)$.
Foata and Zeilberger proved in~\cite{foata4} the following
extension of~(MMT)\,:
$$\text{(FZ)} \qquad
\sum_{\s k = (k_1,\ldots,k_m)} \,\frac{1}{k_1! \cdots k_m!} \,
\sum_{\pi \in \pP(\s k)} \, v_\beta(\pi)
= \left(\frac1{\det(I-A)}\right)^{\beta}.
$$
Note that the right-hand side of~(FZ) is well defined for all complex
values~$\be$, but we will avoid this generalization for simplicity.

Now, in the spirit of Subsection~\ref{ss:basic-noncom} one can
ask whether (FZ) can be extended to a non-commutative setting.
As it turns out, this is quite straightforward given the structure
of our bijection~$\vp$.  As an illustration, we will work in the
setting of Section~\ref{rq}.

\begin{thm}
Let $A = (a_{ij})_{m \times m}$ be a right quantum matrix and
assume that the variables $x_1,\ldots,x_m$ commute with each
other and with all~$a_{ij}$, $1 \le i,j \le m$.
Then, in the above notation, we have:
\begin{equation}\label{e:beta}
\sum_{\s k = (k_1,\ldots,k_m)} \,\frac{1}{k_1! \cdots k_m!} \,
\sum_{\pi \in \pP(\s k)} \, v_\beta(\pi)
= \left(\frac1{\det(I-A)}\right)^{\beta}\,,
\end{equation}
where the summation is over all non-negative integer vectors
$\s k = (k_1,\ldots,k_m)$ and $\det(\cdot)$ is the Cartier-Foata
determinant. \end{thm}

\begin{proof}  We prove the theorem by reduction to
Foata-Zeilberger's identity~(FZ) and our previous results.
First, by Theorem~\ref{rq4}, every term on the right-hand side of
 equation~\eqref{e:beta} is a concatenation of~$\beta$ o-sequences.
 Using bijection~$\vp$ as in the proof of Theorem~\ref{rq4},
 we conclude that
 the sum of all concatenations of~$\beta$ o-sequences is equal
 to a weighted sum of all o-sequences modulo the ideal~$\p I_{\rtq}$.
 In other words, $\left(\det(I-A)\right)^{-\beta}$ \ts is a weighted
 sum of words $v(\pi)/(k_1! \cdots k_m!)$, for $\pi \in \pP(\s k)$,
 and the coefficients are equal to the number of concatenations
 of~$\beta$ o-sequences that are transformed into the given p-sequence.
 Therefore, the coefficients must be the same as in the commutative case.
 Now Foata-Zeilberger's equation~(FZ) immediately implies the theorem.
\end{proof}

\begin{exm} \label{ex:beta}
{\rm
Figure \ref{fig5} illustrates the term
 $(a_{13}a_{22}a_{31})(a_{11}a_{12}a_{23}a_{31})(a_{23}a_{32})$
 in $(\det(I-A))^{-3}$.
 \begin{figure}[ht!]
 \begin{center}
  \includegraphics{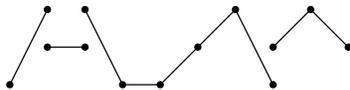}
  \caption{Concatenation of three o-sequences of lengths~3,~4 and~2.}
  \label{fig5}
 \end{center}
 \end{figure}

\noindent
For $\beta=2$, Figure~\ref{fig6} shows  all $(\beta^3+\beta^2)/2=6$ pairs of
 o-sequences whose concatenation  gives the term $a_{11}a_{13}a_{22}a_{31}$
 in $\left(\det(I-A)\right)^{-\beta}$.

 \begin{figure}[ht!]
 \begin{center}
  \includegraphics{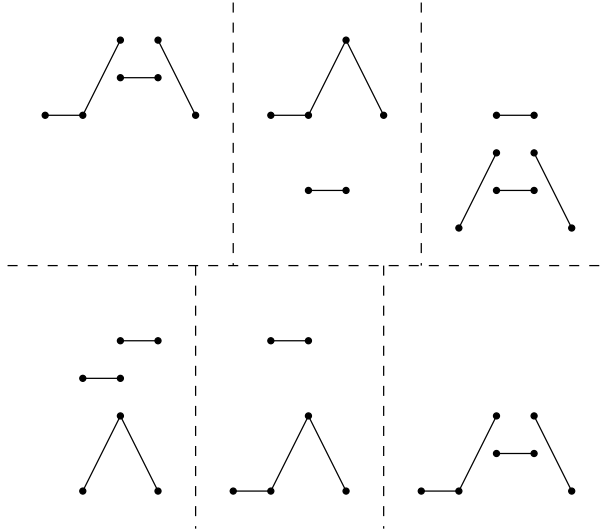}
  \caption{Pairs of o-sequences whose concatenation give
  $a_{11}a_{13}a_{22}a_{31}$ after shuffling.}
  \label{fig6}
 \end{center}
 \end{figure}
}
\end{exm}

\section{Krattenthaler-Schlosser's $q$-analogue} \label{ks}

In the context of multidimensional $q$-series an interesting
$q$-analogue of MacMahon Master Theorem was obtained
in~\cite[Theorem 9.2]{krattenthaler}.  In this section we
place the result in our non-commutative framework and quickly
deduce it from Theorem~\ref{cf1}.

\medskip

We start with some basic definitions and notations.
Let $z_i,b_{ij}$, $1\leq i,j\leq m$, be commutative variables,
and let $q_1,\ldots,q_m \in \C$ be fixed complex numbers.
Denote by~$\p E_i$ the $q_i$-shift operator
$$\p E_i:
\C[z_1,\ldots,z_m] \longrightarrow \C[z_1,\ldots,z_m]
$$
that replaces each occurrence of~$z_i$ by~$q_i \ts z_i$.
We assume that $\p E_r$ commutes with~$b_{ij}$, for all
$1 \le i,j,r \le m$.
For a nonnegative integer vector $\s k = (k_1,\dots,k_m)$,
denote by $[\s z^{\s k}] \ts  F$ the coefficient of
$z_1^{k_1}\cdots z_m^{k_m}$ in the series~$F$.
Denote by~$\s 1$ the constant polynomial~$1$.
Finally, let
$$(a;q)_k \, = \, (1-a) (a-a\ts q) \cdots (1-a\ts q^{k-1}).
$$

\medskip
\begin{thm} [Krattenthaler-Schlosser] \label{t:ks}
Let $A = (a_{ij})_{m\times m}$, where
$$a_{ij} \, = \, z_i\ts \delta_{ij} - z_i \ts b_{ij} \ts \p E_i\,,
 \ \ \ \text{for all} \ \ 1\leq i,j\leq m.
 $$
Then, for non-negative integer vector $\s k =(k_1,\ldots,k_m)$
we have:
\begin{equation} \label{e:ks}
 [\s z^{\s 0}] \, \prod_{i=1}^m \left( \sum_{j=1}^m b_{ij}
 z_j/z_i ; q_i \right)_{k_i} = \left[\s z^{\s k}\right] \, \left(
 \frac1{\det (I -A)} \cdot \s 1 \right),
\end{equation}
 where $\det(\cdot)$ is the Cartier-Foata determinant.
\end{thm}

Note that the right-hand side of~\eqref{e:ks} is non-commutative and
(as stated) does not contain~$q_i$'s, while the left-hand side
contains only commutative variables and~$q_i$'s.
It is not immediately obvious and was shown in~\cite{krattenthaler}
that the theorem reduces to the MacMahon Master Theorem.  Here
we give a new proof of the result.

\medskip
\begin{proof}[Proof of Theorem~\ref{e:ks}]
Think of variables~$z_i$ and~$b_{ij}$ as operators acting on
polynomials by multiplication.  Then a matrix entry~$a_{ij}$
is an operator as well.  Note that multiplication by $z_i$ and
the operator $\p E_j$ commute for $i \neq j$.  This implies that
the equation~\eqref{cf4} holds, i.e.\hspace{-0.07cm} that~$A$ is a Cartier-Foata
matrix. Let $x_1,\ldots,x_m$ be formal variables that commute
with each other and with $a_{ij}$'s.
By Theorem~\ref{cf1}, for the operator on
the right-hand side of~\eqref{e:ks} we have:
$$\frac1{\det (I -A)} \, = \,
\sum_{\s r = (r_1,\ldots,r_m)} \, G(r_1,\ldots,r_m),$$
where
$$G(r_1,\ldots,r_m) \, = \, \left[\s x^{\s r}\right] \,
\prod_{i =\tth 1..{m}}^{\longrightarrow} \,
\bigl(a_{i1}x_1+\ldots+a_{im}x_m\bigr)^{r_i}.$$
Recall that $a_{ij} = z_i \ts(\delta_{ij} - b_{ij}\ts \p E_i)$.
Now observe that every coefficient
$G(r_1,\ldots,r_m) \cdot \s 1$ is equal to
$\s z^{\s r}$ times a polynomial in~$b_{ij}$ and~$q_i$.
Therefore, the right-hand side of~\eqref{e:ks} is equal to
$$\left[{\s z}^{\s k}\right] \, \left(
 \frac1{\det (I -A)} \cdot \s 1 \right) \, = \,
\left[\s z^{\s k}\right] \, \left(\sum_{\s r} \, G(r_1,\ldots,r_m)
\cdot \s 1\right) \, = \,
\left[\s z^{\s k}\right] \, \bigl(  G(k_1,\ldots,k_m)
\cdot \s 1\bigr).
$$
This is, of course, a sum of $[\s z^{\s k}](\alpha \cdot \s 1)$
over all o-sequences $\alpha$ of type $\s k$. Define
$$c_{ij}^k=z_i \delta_{ij} - z_i b_{ij} q_i^{k-1} \ \ \, \text{and}
\ \ d_{ij}^k=z_j \delta_{ij} - z_j b_{ij} q_i^{k-1}.
$$
It is easy to prove by induction that
$$
a_{i\lambda_1} a_{i \lambda_2} \cdots a_{i \lambda_\ell}
\cdot \s 1 \, = \, c_{i \lambda_1}^\ell c_{i \lambda_2}^{\ell - 1}
\cdots c_{i \lambda_\ell}^1.
$$
Therefore, for every o-sequence
\begin{equation} \label{e:o-seq}
 \alpha \, = \, a_{1 \lambda_1^1}a_{1 \lambda_2^1}\cdots
 a_{1 \lambda_{k_1}^1} a_{2 \lambda_1^2}a_{2 \lambda_2^2}\cdots
 a_{2 \lambda_{k_2}^2}  \cdots a_{m \lambda_1^m}a_{m \lambda_2^m}
 \cdots a_{m \lambda_{k_m}^m}
\end{equation}
we have:
$$\aligned
\alpha \cdot \s 1 \, & = \,
c_{1 \lambda_1^1}^{k_1}c_{1 \lambda_2^1}^{k_1-1} \,
\cdots \, c_{1 \lambda_{k_1}^1}^1 \, c_{2 \lambda_1^2}^{k_2}
c_{2\lambda_2^2}^{k_2-1}
\cdots \, c_{2 \lambda_{k_2}^2}^1 \, \cdots \, c_{m \lambda_1^m}^{k_m} \,
c_{m \lambda_2^m}^{k_m-1} \cdots \, c_{m \lambda_{k_m}^m}^1 \\
& = \, d_{1 \lambda_1^1}^{k_1}d_{1 \lambda_2^1}^{k_1-1} \cdots \,
d_{1 \lambda_{k_1}^1}^1 d_{2 \lambda_1^2}^{k_2}
d_{2 \lambda_2^2}^{k_2-1} \cdots \, d_{2 \lambda_{k_2}^2}^1 \cdots \,
d_{m \lambda_1^m}^{k_m}d_{m \lambda_2^m}^{k_m-1} \cdots \,
d_{m \lambda_{k_m}^m}^1,
\endaligned
$$
where the second equality holds because~$\alpha$ is a balanced sequence.
On the other hand,
$$\left[\s z^{\s 0}\right] \,
\prod_{i=1}^m \left( \sum_{j=1}^m b_{ij} z_j/z_i ; q_i \right)_{k_i} \,
= \,
\left[\s z^{\s k}\right] \, \prod_{i=1}^m \,
\prod_{j=1}^{k_i} (d_{i1}^j+\ldots+d_{im}^j)
$$
is equal to the sum of
$$\left[\s z^{\s k}\right]\left(d_{1 \lambda_1^1}^{k_1}
d_{1 \lambda_2^1}^{k_1-1} \cdots \,d_{1 \lambda_{k_1}^1}^1
d_{2 \lambda_1^2}^{k_2}d_{2 \lambda_2^2}^{k_2-1} \cdots \,
d_{2 \lambda_{k_2}^2}^1 \cdots \,d_{m \lambda_1^m}^{k_m}
d_{m \lambda_2^m}^{k_m-1} \cdots \,d_{m \lambda_{k_m}^m}^1\right)$$
over all o-sequences $\alpha$ of form \eqref{e:o-seq}.
This completes the proof. \end{proof}

\section{Proofs of linear algebra propositions} \label{proofs}

\subsection{Proof of Proposition \ref{cf3}}
The proof imitates the standard linear algebra proof
in the commutative case.  We start with the following easy result.
 \begin{lemma}
  Let $B=(b_{ij})_{m \times m}$.
  \begin{enumerate}
   \item If $B$ satisfies \eqref{cf4} and if $B'$ denotes the matrix we
   get by interchanging adjacent columns of $B$, then $\det B'=-\det B$.
   \item If $B$ satisfies \eqref{cf4} and has two columns equal, then $\det B=0$.
   \item If $B^{ij}$ denotes the matrix obtained from $B$ by deleting the
   $i$-th row and the $j$-th column, then
   $$\det B = \sum_{i=1}^m (-1)^{m+i} (\det B^{im}) \tth b_{im}.$$
  \end{enumerate}
 \end{lemma}
 \noindent
The proof of the lemma is completely straightforward.  Now take $B=I-A$ and
recall that $B$ is invertible.  The $j$-th coordinate of the matrix product
 $$
 (\det B^{11},-\det B^{21},\ldots,(-1)^m B^{m1}) \cdot B
 $$
 is $\sum_{i=1}^m (-1)^i \det B^{i1} b_{ij}$.  Since $B$ satisfies \eqref{cf4},
 this is equal to $\det B \cdot \delta_{1j}$ by the lemma. But then
 $$(\det B^{11},-\det B^{21},\ldots,(-1)^m B^{m1}) =
 \det B \cdot (1,0,\ldots,0) \cdot B^{-1}$$
 and
 $$(B^{-1})_{11}=(\det B)^{-1} \cdot \det B^{11}. \eqno \qed
 $$

\subsection{Proof of Proposition \ref{rq3}}
Follow the same scheme as in the previous subsection.
The following is a well-known result
(see e.g. \cite[Lemmas 2.3 and 2.4]{garoufalidis}
or \cite[Properties 5 and 6]{foata2}).

 \begin{lemma}
  Let $B=(b_{ij})_{m \times m}$.
  \begin{enumerate}
   \item If $B$ satisfies \eqref{rq2} and if $B'$ denotes the matrix
   we get by interchanging adjacent columns of $B$, then $\det B'=-\det B$.
   \item If $B$ satisfies \eqref{rq2} and has two columns equal, then $\det B=0$.
   \item If $B^{ij}$ denotes the matrix obtained from $B$ by deleting the
   $i$-th row and the $j$-th column, then
   $$\det B = \sum_{i=1}^m (-1)^{m+i} (\det B^{im})\tth b_{im}. \eqno \qed$$
  \end{enumerate}
 \end{lemma}
 The rest follows verbatim the previous argument.

\subsection{Proof of Proposition \ref{q6}}

Foata and Han introduced (\cite[Section 3]{foata1}) the so-called ``$1=q$
principle'' to derive identities in the algebra $\p A/\p I_{q-\rtq}$ from
those in the algebra $\p A/\p I_{\rtq}$.

\begin{lemma} \label{q7} {\rm (``$1=q$ principle'')}
 Let $\phi \colon \p A \to \p A$ denote the linear map induced by
 $$\phi\left(a_{\lambda,\mu}\right) =
 q^{\inv \mu - \inv \lambda} a_{\lambda,\mu}.
 $$
 Then:
 \begin{enumerate}
  \renewcommand{\labelenumi}{(\alph{enumi})}
  \item $\phi$ maps $\p I_{\rtq}$ into $\p I_{q-\rtq}$
  \item Call $a_{\lambda,\mu}$ a \emph{circuit} if $\lambda$ is a rearrangement
  of $\mu$ (i.e.\hspace{-0.07cm} if $\lambda$ and $\mu$ contain the same letters with the
  same multiplicities). Then $\phi(\alpha \beta)=\phi(\alpha)\phi(\beta)$
  for $\alpha,\beta$ linear combinations of circuits.
 \end{enumerate}
\end{lemma}

We include the proof of the lemma since we need to generalize it later on.

\medskip

\begin{prf}
 (a) \, It suffices to prove the claim for elements of the form
 $$\alpha=a_{\lambda,\mu}(a_{ik}a_{jk} - a_{jk}a_{ik})a_{\lambda',\mu'}$$
 and
 $$
 \beta=a_{\lambda,\mu}(a_{ik}a_{jl}-a_{jk}a_{il}
 - a_{jl}a_{ik} + a_{il}a_{jk})a_{\lambda',\mu'}
 $$
 with $i < j$ (and $k < l$). Note that the sets of inversions of the
 words $\lambda i j \lambda'$ and $\lambda i j \lambda'$ differ only in the
 inversion $(i,j)$. Therefore $\phi(\alpha)$ is a multiple
 of
 $$a_{ik}a_{jk} - q^{-1} a_{jk}a_{ik}.$$
 For the~$\beta$ the proof is analogous.

\smallskip

 (b) \, It suffices to prove the claim
 for $\alpha,\beta$ circuits, i.e.\hspace{-0.07cm}
 $\alpha=a_{\lambda,\mu}$ with $\lambda$ a rearrangement of $\mu$ and
 $\beta=a_{\lambda',\mu'}$ with $\lambda'$ a rearrangement of $\mu'$.
 The set of inversions of $\lambda \lambda'$ consists of the inversions
 of $\lambda$, the inversions of $\lambda'$, and the pairs $(i,j)$ where $\lambda_i>\mu_j$.
 Similarly, the set of inversions of $\mu \mu'$ consists of the inversions
 of $\mu$, the inversions of $\mu'$, and the pairs $(i,j)$ where $\lambda'_i>\mu'_j$. Since $\lambda$
 is a rearrangement of $\mu$ and $\lambda'$ is a rearrangement of $\mu'$,
 $\inv(\mu \mu') - \inv (\lambda \lambda')=
 (\inv \mu - \inv \lambda) + (\inv \mu' - \inv \lambda')$,
 which concludes the proof.\qed
\end{prf}

By Proposition~\ref{rq3}, we have:
$$\det(I-A) \cdot \left((I-A)^{-1}\right)_{11} -
\det\left(I-A^{11}\right) \in \p I_{\rtq}\,.
$$
It is clear that
$$\phi(\det (I-A))\!=\!\phi\left( \sum (-1)^{|J|} \det A_J \right)\!
= \!\sum (-1)^{|J|} \detq A_J \!=\!\detq (I-A),
$$
where the sums go over all subsets $J \subseteq \set{1,\ldots,m}$.
Similarly,
$$\phi\left( \left((I-A)^{-1}\right)_{11} \right) =
\left((I-\widetilde A)^{-1}\right)_{11}.
$$
Now the result follows from Lemma~\ref{q7}.\qed

\subsection{Proofs of Propositions \ref{cfq7}, \ref{cfqij8} and \ref{qij7}}
 \label{ss:proofs-qij}
 The result can be derived from Propositions \ref{cf3} and \ref{rq3} by a simple
 extension of the ``$1=q$ principle''.

 \begin{lemma} \label{proofs1} {\rm (``$1=q_{ij}$ principle'')}
  Call an element $\sum_{i \in \p J} c_i a_{\lambda_i,\mu_i}$ of $\p A$
  \emph{balanced} if for any $i,j \in \p J$, $\lambda_i$ is a reshuffle
  of~$\lambda_j$ and $\mu_i$ is a reshuffle of $\mu_j$.\\
  Let $\phi \colon \p A \to \p A$ denote the linear map induced by
  $$\phi\left(a_{\lambda,\mu}\right) =
  \left(\prod_{(i,j) \in I(\mu)} q_{\mu_j \mu_i} \prod_{(i,j) \in I(\lambda)}
  q_{\lambda_j \lambda_i}^{-1}\right) a_{\lambda,\mu}.
  $$
  Choose a set $\p S$ with balanced elements, denote by $\p I$ the ideal
  generated by $\p S$, and by $\p I_{\s q}$ the ideal generated by $\phi(\p S)$.
  Then
  \begin{enumerate}
   \renewcommand{\labelenumi}{(\alph{enumi})}
   \item $\phi$ maps $\p I$ into $\p I_{\s q}$,
   \item $\phi(\alpha \beta)=\phi(\alpha)\phi(\beta)$ for $\alpha,\beta$
   linear combinations of circuits.
  \end{enumerate}
 \end{lemma}

\noindent
 The proof of lemma follows verbatim the proof of Lemma~\ref{q7}.
Propositions \ref{cfq7} and \ref{cfqij8} follow from Proposition \ref{cf3},
and Proposition \ref{qij7} follows from Proposition \ref{rq3}.
We omit the details.  \qed

\section{Final remarks} \label{s:final}

\subsection{} \label{ss:final-koszul}
A connection between Cartier-Foata free partially-commutative monoids
and Koszul duality was established by Kobayashi~\cite{kob} and can
be stated as follows.  Let~$G$ be a graph on~$[n] = \{1,\ldots,n\}$.
Consider a quadratic algebra $\p A_G$ over~$\C$ with variables
$x_1,\ldots,x_n$ and relations $x_ix_j = x_j x_i$ for every
edge~$(i,j) \in G$, $i \ne j$, and $x_i^2 = i$ if there is a
loop at~$i$.  It was shown by Fr\"oberg in full generality
that~$\p A_G$ is Koszul, and the Koszul dual algebra $\p A^!_G$ has
a related combinatorial structure (see~\cite{Fr}).  This generalizes the
classical case of a complete graph $G=K_n$, where $\p A_G$ is a symmetric
and $\p A^!_G$ is an exterior algebra.  We refer to~\cite{pol} for a
extensive recent survey on quadratic algebras and Koszul duality.

Now, Kobayashi observed that one can view the Cartier-Foata M\"obius
inversion theorem for the partially commutative monoid corresponding
to a graph~$G$ (see~\cite{cf}) as a statement about Hilbert series:
\begin{equation} \label{final1}
 A_G(t) \cdot A_G^!(t) = 1
\end{equation}
where $A(t) = \sum_i \dim \p A^i t^i $ for a graded algebra $\p A = \oplus \p A^i$.
In effect, Kobayashi gives an explicit construction of the Koszul
complex for~$\p A_G$ by using Cartier-Foata's involution~\cite{kob}.

Most recently, Hai and Lorenz made a related observation, by showing
that one can view the Master Theorem as an identity of the same type
as \eqref{final1} but for the characters rather than dimensions~\cite{hai}.
This allowed them to give an algebraic proof of the
Garoufalidis-L\^e-Zeilberger theorem.  In fact, they present
a general framework to obtain versions of the Master Theorem for
other Koszul algerbras (which are necessarily quadratic) and a
(quantum) group acting on it.

\subsection{} \label{ss:final-manin}  From our presentation, one
may assume that the choice of a $(q_{ij})$-analogue was a lucky
guess or a carefully chosen deformation designed to make the technical
lemmas work.  This was not our motivation, of course.  These
quadratic algebras are well known generalizations of the classical
quantum groups of type~$A$ (see~\cite{manin1,manin2,manin-book}).
They were introduced and extensively studied by Manin, who also
proved their Koszulity and defined the corresponding (generalized)
quantum determinants.

While our proof is combinatorial, we are confident that the
Hai-Lorenz approach will work in the $(q_{ij})$-case
as well.  While we do not plan to further investigate this
connection, we hope the reader find it of interest to pursue
that direction.

\subsection{} \label{ss:final-quasi}  For matrices over general
rings, the elements of the inverse matrices are called
\emph{quasi-determinants}~\cite{gelfand} (see also~\cite{gelfand2}).
They were
introduced by Gelfand and Retakh, who showed that in various
special cases these quasi-determinants are the ratios of two
(generalized) determinants.  In particular, in the context
of non-commu\-ta\-tive determinants they established
Propositions~\ref{cf3}, \ref{q6} and a (slightly weaker)
corresponding result for the super-analogue.

In a more general context, Etingof and Retakh showed the
analogue of this result for all twisted quantum groups~\cite{ER}.
Although they do not explicitly say so, we believe one can probably
deduce our most general Proposition~\ref{qij7} from~\cite{ER}
and the above mentioned Manin's papers.  Interestingly, it follows
from~\cite{ER} and our work that the (non-commutative) determinants
of minors considered in this paper always commute with each other.
We do not need this observation for our telescoping argument.

Let us mention here that the inverse matrix $(I-A)^{-1}$
appears in the same context as in this paper in the study
of quasi-determinants~\cite{gelfand2} and the
non-commutative Lagrange inversion~\cite{pak}.

\subsection{} \label{ss:final-super}  The relations for variables
in our super-analogue are somewhat different from those studied
in the literature (see e.g.~\cite{manin-book}).  Note also that
our super-determinant is different from the
\emph{Berezinian}~\cite{Ber} (see also~\cite{gelfand2,manin1}).
We are somewhat puzzled by this and hope to obtain the ``real''
super-analogue in the future.

\subsection{} \label{ss:final-etingof}  The relations
 studied in this paper always lead to quadratic algebras.
While the deep reason lies in the Koszul duality, the fact
that Koszulity can be extended to non-quadratic algebras
is suggestive~\cite{Berger}.  The first such effort is made
in~\cite{EP} where an unusual algebraic extension of
MacMahon Master Theorem is obtained.

\subsection{} \label{ss:final-beta}  While we do not state
the most general result combining both $\beta$-extension and
$(q_{ij})$-analogue, both the statement and the proof follow
verbatim the presentation in Section~\ref{beta}.  Similarly,
the results easily extend to all complex values~$\be \in \C$.

Let us mention
here that the original $\beta$-extension of the Master Theorem
(given in~\cite{foata4}) follows easily from the $\beta$-extension
of the Lagrange inversion~\cite{zeng}.  In fact, the proof of the
latter is bijective.

\subsection{} \label{ss:final-bozon} In the previous
papers~\cite{foata1,foata2,foata3,garoufalidis} the authors used
$Boz(\cdot)$ and $Fer(\cdot)$ notation for the left- and the right-hand side
of~\eqref{e:macm-D}.  While the implied connection is not unjustified,
it might be misleading when the results are generalized.  Indeed, in
view of Koszul duality connection (see Subsection~\ref{ss:final-koszul}
above) the algebras can be interchanged, while giving the same result
with notions of \emph{Bozon} and \emph{Fermion} summations switched.
On the other hand, we should point out that in the most interesting
cases the \emph{Fermion} summation is finite, which makes it special
from combinatorial point of view.

\subsection{} \label{ss:kratt-schloss}
The Krattenthaler-Schlosser's $q$-analogue (Theorem~\ref{t:ks}) is
essentially a byproduct of the author's work on $q$-series.  It was
pointed out to us by Michael Schlosser that the Cartier-Foata
matrices routinely appear in the context of ``matrix inversions''
for $q$-series (see~\cite{krattenthaler,Sc}).  It would be
interesting to see if our extensions
(such as Cartier-Foata $q_{ij}$-case in Section~\ref{cfqij}) can
can be used to obtain new results, or give new proofs of existing
results.

\medskip

\subsection*{Acknowledgements}

\noindent
We are grateful to Alexander Polishchuk, Vic Reiner, Vladimir Retakh,
Michael Schlosser, Richard Stanley and Doron Zeilberger for the
interesting discussions and help with the references.  We are especially
thankful to Pavel Etingof and Christian Krattenthaler for helping
us understand the nature of Theorem~\ref{t:ks}. The second author
was partially supported by the NSF.


\end{document}